\documentclass[titlepage,twoside,12pt]{article}
\usepackage{amssymb}
\usepackage{amsfonts}
\begin{document}

{\LARGE \bf Two Essays on the Archimedean \\ \\ versus Non-Archimedean Debate} \\ \\

{\bf Elem\'{e}r E ~Rosinger} \\ \\
{\small \it Department of Mathematics \\ and Applied Mathematics} \\
{\small \it University of Pretoria} \\
{\small \it Pretoria} \\
{\small \it 0002 South Africa} \\
{\small \it eerosinger@hotmail.com} \\ \\

{\LARGE \bf Part 1} \\ \\

{\LARGE \bf Is Space-Time Non-Archimedean ?} \\ \\

{\bf Abstract} \\

For more than two millennia, ever since Euclid's geometry, the so called Archimedean Axiom has been accepted without sufficiently explicit awareness of
that fact. The effect has been a severe restriction of our views of space-time, a restriction which above all affects Physics. Here it is argued that,
ever since the invention of Calculus by Newton, we may actually have {\it empirical evidence} that time, and thus space as well, are not Archimedean. \\ \\

{\bf 1. A Brief Review of the Axioms of Euclidean Geometry} \\

Ever since the discoveries of non-Euclidean geometries by Lobachevski and Bolyai in the early 1800s, the axioms of Euclidean geometry have been divided in
two : on one hand, one has the Axiom of Parallels, while the rest of the axioms constitutes what is called Absolute Geometry. And the non-Euclidean
geometries occur when the axioms of Absolute Geometry are augmented either by the axiom of inexistence of parallel lines, or by the axiom of the
existence of more than one parallel line. \\

What happens, however, is that the axioms of Absolute Geometry contain what came to be called the Archimedean Axiom, which in its simplest formulation
for the real line $\mathbb{R}$ is as follows \\

(1.1)~~~ $ \exists~~ u > 0 ~:~ \forall~~ x > 0 ~:~ \exists~~ n \in \mathbb{N} ~:~ n u > x $ \\

In other words, there exists a "step-length" $u > 0$, so that, no matter how far away a point $x > 0$ would be, one can starting at the point $0$ go
beyond it in a finite number $n$ of "steps". Of course, in the case of the real line $\mathbb{R}$, every strictly positive $u > 0$ is such a
"step-length". The issue with the Archimedean Axiom, however, is that there should exist at least one such "step-length". \\

Let us therefore denote by AG the axioms of Absolute Geometry, including thus the Archimedean Axiom, while by PAG, that is, Pure Absolute Geometry, we shall denote the axioms of Absolute Geometry {\it minus} the Archimedean Axiom. In this way, we have \\

(1.2)~~~ $ AG = PAG + AA $ \\

where by AA we denote the Archimedean Axiom. \\

In the sequel we shall be interested in PAG. Furthermore, we shall argue that Calculus, more precisely, the mental process of implementing Calculus can offer an empirical evidence for the fact that space-time is rather described by PAG, than by AG. \\ \\

{\bf 2. Happily and for Evermore Wallowing in the \\
        \hspace*{0.4cm} "Egyptian Bondage" of the Archimedean Axiom ?} \\

Needless to say, the utility of the Archimedean Axiom had been vital practically in, among others, the perennial redrawing of farm land in ancient Egypt, after the yearly floods of the Nile. \\

Yet its uncritical, and in fact, less than fully conscious, nevertheless, longtime acceptance has incredibly limited our perception and intuition of
space-time. \\
And here we are even nowadays, with such a most severe limitation, and we are not aware of the fact that our liberation from it does not require the
additional acceptance of anything at all, but on the contrary, the mere setting aside of an assumption, an assumption which somehow got stuck to us in
ancient times, an assumption which ever since we took for granted as being a cornerstone of the way space-time is, namely, the Archimedean Axiom. \\

And as with all such deep seated assumptions, it is not at all easy to follow the call for its setting aside. And it is not easy, even if the consequent immense opening in the wealth of possibilities for space-time does not require the acceptance of any new assumption, but only the setting aside of an old one. \\

And if it is not easy for mathematicians to follow such an immense opening, one need not wonder that physicists may not find it easier to do so ... \\

Yet rather simple fundamental issues relating to space-time, [2], keep now and then coming to the surface as if in a somewhat surprising manner ... \\ \\

{\bf 3. But Now, What Can Non-Archimedean Space-Time Do \\
        \hspace*{0.4cm} for You ?} \\

As it happens, more than three centuries earlier, and more than a century prior to the emergence of non-Euclidean geometries, Leibniz introduced
"infinitesimals" in the structure of the real line $\mathbb{R}$. Unfortunately, there was not enough Mathematics at the time to do that in a systematic
and rigorous enough manner. However, several of its uses in Calculus as suggested by Leibniz proved to be not only correct, but remarkably efficient in
their brevity, as well as in their immediate intuitive clarity. The idea of "infinitesimals" was in fact a consequence of one of the foundational
philosophical ideas of Leibniz, namely, of "monads". And in this regard, the successful use he could make of "infinitesimals" in Calculus must have quite
impressed him as one likely confirmation of the reality and foundational nature of "monads". \\

The systematic, rigorous, and also surprisingly far reaching mathematical formulation of "infinitesimals" was given by Abraham Robinson's 1966
construction of Nonstandard Analysis which introduced, among others, time and space structures significantly more rich and complex than the classical
ones, see section 4 below. \\
The fact that Nonstandard Analysis has nevertheless failed to become popular even within the mathematical community is due to a cost-return situation in
which the vast majority of mathematicians have, rightly or wrongly, decided that the returns do in no way justify the costs, namely, what to so many
appear to be excessive technical complication. \\

{\bf Reduced Power Algebras, or RPA-s} \\

Fortunately however, as it happens, the very same enriched time and space structures given by Nostandard Analysis - and called {\it reduced power
algebras}, or in short, RPA-s - can be obtained in a far simpler manner, namely, using only 101 Algebra, and specifically, the concepts of ring, ideal,
quotient, [1], as well as the rather simple and intuitive set theoretic concept of filter. \\
What one loses by that much more simple, easy and direct approach is the Transfer Principle in Nonstandard Analysis. Furthermore, as is well known, the
Transfer Principle in Nonstandard Analysis suffers itself from severe limitations, namely, it is restricted to entities which can be described by what is
called First Order Predicate Logic. On the other hand, a large amount of entities in Calculus, not to mention the rest of Mathematics, cannot be
formulated within First Order Predicate Logic, thus fall outside of the range of applicability of the Transfer Principle. \\

However, regarding one's full liberation from the "Egyptian Bondage", that is, entering in the realms of PAG, such a loss does {\it not} inconvenience to
any significant extent. In this way, the RPA-s offer remarkable instances of PAG, thus we can see the situation as given by the inclusion \\

(3.1)~~~ $ RPA ~\subseteq~ PAG $ \\

{\bf Constructing RPA-s} \\

We recall here briefly the construction of the RPA-s which have been well known in Model Theory, where they constitute one of the most important basic
concepts. \\ Notations and details used in the sequel can be found in [4-6]. \\

This construction happens in the following {\it three} steps

\begin{itemize}

\item first, one chooses an arbitrary {\it infinite} index set $\Lambda$ and constructs the power algebra $\mathbb{R}^\Lambda$ which is but the set of
all functions $f : \Lambda \longrightarrow \mathbb{R}$ considered with the usual addition and multiplication of functions,

\item second, one chooses a {\it proper} ideal ${\cal I}$ in the power algebra $\mathbb{R}^\Lambda$,

\item third, one constructs the quotient algebra

\end{itemize}

(3.1)~~~ $ \mathbb{A} = \mathbb{R}^\Lambda / {\cal I} $ \\

which is called the {\it reduced power algebra}, or in short, RPA. \\

An important simplification of this construction can be obtained by replacing proper ideals ${\cal I}$ in the algebra $\mathbb{R}^\Lambda$, with the
simpler mathematical structures of {\it filters} ${\cal F}$ on $\Lambda$. Here we recall that a {\it filter} ${\cal F}$ on $\Lambda$ is a set of subsets
$I \subseteq \Lambda$ with the following four properties \\

(3.2)~~~ $ {\cal F} \neq \phi $ \\

(3.3)~~~ $ \phi \notin {\cal F} $ \\

(3.4)~~~ $ I, J \in {\cal F} \Longrightarrow I \cap J \in {\cal F} $ \\

(3.5)~~~ $ I \in {\cal F},~~ I \subseteq J \subseteq \Lambda \Longrightarrow J \in {\cal F} $ \\

Thus such filters can be seen as collections of {\it large} subsets of $\Lambda$. Indeed, (3.2) means that there exist such large subsets, and certainly,
none of them is void, as required by (3.3). Condition (3.4) means that the intersection of two large subsets is still a large subset, while
(3.5) simply means that a subset which contains a large subset is itself large. In particular, $\Lambda$ itself is large, thus $\Lambda \in {\cal F}$. \\

And important example of filter on $\Lambda$ is the {\it Frech\'{e}t} filter, given by \\

(3.6)~~~ $ {\cal F}_{re} ( \Lambda ) = \{~ I \subseteq \Lambda ~~|~~ \Lambda \setminus I ~~\mbox{is finite} ~\} $ \\

The mentioned simplification comes about the following one-to-one simple correspondence between proper ideals ${\cal I}$ in the algebra
$\mathbb{R}^\Lambda$ and filters ${\cal F}$ on $\Lambda$, namely \\

(3.7)~~~ $ {\cal I} ~~\longmapsto~~ {\cal F}_{\cal I} = \{ Z ( x ) ~|~ x \in {\cal I} \} ~~\longmapsto~~ {\cal I}_{{\cal F}_{\cal I}} = {\cal I} $ \\

(3.8)~~~ $ {\cal F} ~~\longmapsto~~ {\cal I}_{\cal F} = \{ x \in \mathbb{R}^\Lambda ~|~ Z ( x ) \in {\cal F} \}
                         ~~\longmapsto~~ {\cal F}_{{\cal I}_{\cal F}} = {\cal F}$ \\

where for $x \in \mathbb{R}^\Lambda$ we denote $Z ( x ) = \{ \lambda \in \Lambda ~|~ x ( \lambda ) = 0 \}$. \\

{\bf RPA-s Extend the Real Line $\mathbb{R}$} \\

An important consequence of (3.1), (3.7), (3.8) is that the mapping \\

(3.9)~~~ $ \mathbb{R} \ni r \longmapsto u_r + {\cal I} \in A = \mathbb{R}^\Lambda / {\cal I} $ \\

is an {\it injective algebra homomorphism} for every proper ideal ${\cal I}$ in $\mathbb{R}^\Lambda$, where $u_r \in \mathbb{R}^\Lambda$ is defined by
$u_r ( \lambda ) = r$, for $\lambda \in \Lambda$. Indeed, in view of (3.7), (3.8), we have for $r \in \mathbb{R}$ \\

(3.10)~~~ $ u_r \in {\cal I} ~~\Longrightarrow~~ r = 0 $ \\

since ${\cal I} = {\cal I}_{{\cal F}_{\cal I}}$, while $u_r \in {\cal I}_{{\cal F}_{\cal I}}$ gives $Z ( u_r ) \in {\cal F}_{\cal I}$, thus $Z ( u_r )
\neq \phi$, which means $r = 0$. \\

Also (3.7), (3.8) imply that the RPA-s (3.1) can be represented in the form \\

(3.11)~~~ $ \mathbb{A}_{\cal F} = \mathbb{R}^\Lambda / {\cal I}_{\cal F} $ \\

where ${\cal F}$ ranges over all the filters on $\Lambda$, thus the {\it injective algebra homomorphisms} (3.9) become \\

(3.12)~~~ $ \mathbb{R} \ni r \longmapsto u_r + {\cal I}_{\cal F} \in \mathbb{A}_{\cal F} $ \\

which means that the usual real line $\mathbb{R}$, which is a field, is in fact a {\it subalgebra} of each of the RPA-s, see (3.1), (3.9), (3.11), (3.12),
namely \\

(3.13)~~~ $ \mathbb{R} ~\subseteq~ \mathbb{A}_{\cal F} $ \\

However, we have to note that for some filters ${\cal F}$ on $\Lambda$, one may have equality in (3.13), which of course, is of no interest as far as
RPA-s are concerned. Another trivial case to avoid is that of filters generated by a nonovid set $I \subseteq \Lambda$, namely, ${\cal F}_I = \{~ J
\subseteq \Lambda ~~|~~ J \supseteq I ~\}$, when we obtain \\

(3.14)~~~ $ \mathbb{A}_{{\cal F}_I} = \mathbb{R}^I $ \\

thus instead of a reduced power algebra, we only have a power algebra. In this regard, it is easy to show that, in (3.13) we have the strict inclusion \\

(3.15)~~~ $ \mathbb{R} ~\subsetneqq~ \mathbb{A}_{\cal F} $ \\

and also avoid the case of (3.14), if and only if \\

(3.16)~~~ $ {\cal F} \supseteq {\cal F}_{re} ( \Lambda ) $ \\

Consequently, from now on, we shall assume that all filters on $\Lambda$ satisfy condition (3.16). \\

{\bf Reduced Power Fields, or RPF-s} \\

Now we turn to {\it reduced power fields}, or in short RPF-s, we are particular cases of RPA-s. Of interest in this respect are a particular case of
filters on $\Lambda$, called {\it ultrafilters} ${\cal U}$, and which are characterized by the condition \\

(3.17)~~~ $ \forall~~ I \subseteq \Lambda ~:~ I \notin {\cal U} \Longrightarrow \Lambda \setminus I \in {\cal U} $ \\

One of their properties relevant in the sequel is that, through (3.7), (3.8), ultrafilters are in one-to-one correspondence with {\it maximal} ideals in
$\mathbb{R}^\Lambda$, namely \\

(3.18)~~~ $ {\cal F} ~\mbox{ultrafilter} ~~\Longrightarrow~~ {\cal I}_{\cal F} ~\mbox{maximal ideal} $ \\

(3.19)~~~ $ {\cal I} ~\mbox{maximal ideal} ~~\Longrightarrow~~ {\cal F}_{\cal I} ~\mbox{ultrafilter} $ \\

For our purposes, it is useful to distinguish between {\it fixed}, and on the other hand, {\it free} ultrafilters on $\Lambda$. The fixed ones are of the
form ${\cal U}_\lambda = \{ I \subseteq \Lambda ~|~ \lambda \in I \}$, for a given fixed $\lambda \in \Lambda$, while the free ones are all the other
ultrafilters on $\Lambda$. It is easy to see that an ultrafilter ${\cal U}$ on $\Lambda$ if free, if and only if \\

(3.20)~~~ $ {\cal F}_{re} ( \Lambda ) \subseteq {\cal U} $ \\

and the existence of free ultrafilters results from the Axiom of Choice. \\

Now we recall from Algebra that \\

(3.21)~~~ $ {\cal I} ~\mbox{maximal ideal in}~ \mathbb{R}^\Lambda ~~\Longleftrightarrow~~ \mathbb{R}^\Lambda / {\cal I} ~\mbox{field} $ \\

And then (3.11), (3.18) - (3.21) will result in \\

(3.22)~~~ $ \mathbb{A}_{\cal F} ~\mbox{field} ~~\Longleftrightarrow~~ {\cal F} ~\mbox{ultrafilter} $ \\

However, in case of a fixed ultrafilter ${\cal F} = {\cal U}_\lambda = \{ I \subseteq \Lambda ~|~ \lambda \in I \}$, for a given fixed $\lambda \in
\Lambda$, it follows that, see (3.14) \\

(3.23)~~~ $ \mathbb{A}_{\cal F} = \mathbb{R} $ \\

thus the corresponding reduced power fields do not lead beyond the field of the usual real line $\mathbb{R}$. Consequently, in (3.22), we shall only be
interested in {\it free ultrafilters} ${\cal F}$ on $\Lambda$, see (3.15), (3.16), (3.20), in which case (3.13) becomes a {\it strict inclusion} of
fields, namely \\

(3.24)~~~ $ \mathbb{R} \subsetneqq \mathbb{A}_{\cal F} $ \\

In the sequel, when dealing with RPF-s, we shall always assume that they correspond to free ultrafilters, thus (3.20) holds for them. \\

{\bf Partial Order on RPA-s} \\

Given, see (3.11), $\mathbb{A}_{\cal F} = \mathbb{R}^\Lambda / {\cal I}_{\cal F}$ a RPA, one can define on it a {\it partial order} $\leq$ by \\

(3.25)~~~ $ \xi ~\leq~ \eta ~~\Longleftrightarrow~~
                \{~ \lambda \in \Lambda ~~|~~ x ( \lambda ) \leq y ( \lambda ) ~\} \in {\cal F} $ \\

where $\xi = x + {\cal I}_{\cal F},~ \eta = y + {\cal I}_{\cal F} \in \mathbb{A}_{\cal F}$, for suitable $x, y \in \mathbb{R}^\Lambda$. The important
fact is that this partial order is compatible with the algebra structure of $\mathbb{A}_{\cal F}$. In particular \\

(3.26)~~~ $ \xi ~\leq~ \eta ~~\Longrightarrow~~ \xi + \chi ~\leq~ \eta + \chi $ \\

for $\xi, \eta, \chi \in \mathbb{A}_{\cal F}$. Also \\

(3.27)~~~ $ \xi ~\leq~ \eta ~~\Longrightarrow~~ \chi \xi ~\leq~ \chi \eta $ \\

for $\xi, \eta, \chi \in \mathbb{A}_{\cal F},~ \chi \geq 0$. Furthermore, the partial order $\leq$, when restricted to $\mathbb{R}$, see (3.13), (3.15),
conicides with the usual total order on $\mathbb{R}$. \\

{\bf Total Order on RPF-s} \\

In the particular case of RPF-s, with their corresponding free ultrafilters, it turns out that the partial order (3.22) is in fact a {\it total order},
namely, for every $\xi, \eta \in \mathbb{A}_{\cal F}$, we have \\

(3.28)~~~ either~ $ \xi \leq \eta $, ~or~ $ \eta \leq \xi $ \\

{\bf Infinitesimals in RPA-s} \\

An element $\xi \in \mathbb{A}_{\cal F}$ is called {\it infinitesimal}, if and only if, for every $r \in \mathbb{R},~ r > 0$, we have in the sense of
(3.25) \\

(3.29)~~~ $ - r \leq \xi \leq r $ \\

Obviously, $\xi = 0 \in \mathbb{A}_{\cal F}$ is trivially an infinitesimal. And in $\mathbb{R}$, the only infinitesimal is $0$. However, the existence of
nonzero infinitesimals in $\mathbb{A}_{\cal F}$ can easily be proved. The set of all infinitesimals, including $0$, in a given $\mathbb{A}_{\cal F}$ is
denoted by \\

(3.30)~~~ $ monad_{\cal F} ( 0 ) $ \\

and recalling Leibniz, it is called the {\it monad} at $0$ in  $\mathbb{A}_{\cal F}$. Usually, when there is no confusion, the index ${\cal F}$ will be
dropped in (3.30). \\

An essential feature of RPA-s is precisely the presence of infinitesimals. \\

{\bf Finite Elements in RPA-s} \\

An element $\xi \in \mathbb{A}_{\cal F}$ is called {\it finite}, if and only if there exists $r \in \mathbb{R},~ r > 0$, such that \\

(3.31)~~~ $ - r \leq \xi \leq r $ \\

The set of all such elements is denoted by \\

(3.32)~~~ $ Gal_{\cal F} ( 0 ) $ \\

and is called the {\it Galaxy} at $0$. As with the monads, when there is no confusion, we shall drop the index ${\cal F}$. \\

It follows easily that \\

(3.33)~~~ $ Gal_{\cal F} ( 0 ) = \bigcup_{r \in \mathbb{R}} monad ( r ) = \mathbb{R} + monad ( 0 ) $ \\

where we denoted \\

(3.34)~~~ $ monad ( r ) = r + monad ( 0 ) $ \\

{\bf Infinitely Large Elements in RPA-s} \\

A consequence of the existence of non-zero infinitesimals in RPA-s, is the presence of {\it infinitely large} elements $\xi \in \mathbb{A}_{\cal F}$,
defined by the condition \\

(3.35)~~~ $ \xi \leq - r $ ~or~ $ r \leq \xi $ \\

for all $r \in \mathbb{R},~ r > 0$. Obviously, the set of all infinitely large elements is given by \\

(3.36)~~~ $ \mathbb{A}_{\cal F} ~\setminus~ Gal_{\cal F} ( 0 ) $ \\

A further essential feature of RPA-s is the presence of infinitely large elements with which one can rigorously perform all the operations in an algebra,
that is, unrestricted addition, subtraction, multiplication and division. With respect to division, we obviously have the {\it injective order reversing}
mapping \\

(3.37)~~~ $ \mathbb{A}_{\cal F} ~\setminus~ Gal ( 0 ) \ni t
                     \longmapsto \frac{1}{t} \in monad ( 0 ) \setminus \{ 0 \} $ \\

As seen in section 4 below, in the particular case of RPF-s, the mapping (3.37) is in fact {\it bijective}, thus establishing one of the many nontrivial
{\it self-similarities} which give the rich and complex structure of RPF-s. \\

{\bf RPA-s Algebras are Non-Archimedean} \\

The presence of nonzero infinitesimals, and thus of infinitely large elements in RPA-s implies that such algebras must be
non-Archimedean. \\

Here, for simplicity, we shall illustrate that implication in the particular case of RPF-s which, as seen,
are totally ordered. And for further convenience, we shall take the index set $\Lambda = \mathbb{N}$. Let therefore ${\cal
U}$ be a free ultrafilter on $\mathbb{N}$ and consider the corresponding RPF \\

(3.38)~~~ $ \mathbb{F}_{\cal U} = \mathbb{R}^\mathbb{N} / {\cal I}_{\cal U} $ \\

We prove now that condition (1.1) which defines in this case the Archimedean property does {\it not} hold. Indeed, let us
assume that it holds, then in view of (3.38), we have \\

(3.39)~~~ $ u = ( u_1, u_2, u_3, \ldots ) + {\cal I}_{\cal U} \in \mathbb{F}_{\cal U} $ \\

where $u_1, u_2, u_3, \ldots \in \mathbb{R}$. If we take now any \\

(3.40)~~~ $ x = ( x_1, x_2, x_3, \ldots ) + {\cal I}_{\cal U} \in \mathbb{F}_{\cal U} $ \\

where $x_1, x_2, x_3, \ldots \in \mathbb{R}$, then (1.1) gives $n \in \mathbb{N}$, such that, see (3.25), $x \leq n u$,
that is \\

(3.41)~~~ $ \{~ \nu \in \mathbb{N} ~~|~~ x_\nu \leq n u_\nu ~\} \in {\cal U} $ \\

However, since ${\cal U}$ is a free ultrafilter, it follows that the set $I = \{~ \nu \in \mathbb{N} ~~|~~ x_\nu \leq
n u_\nu ~\}$ is infinite. And then, we obtain a contradiction, since $x_1, x_2, x_3, \ldots \in \mathbb{R}$ can be
arbitrary, thus in particular, we can choose $x$ so that we have \\

(3.42)~~~ $ \sup_{\nu \in I} \frac{x_\nu}{u_\nu} = \infty $ \\ \\

{\bf 4. The Rich and Complex Self-Similar Structure of \\
        \hspace*{0.5cm} Reduced Power Fields} \\

We recall that the set $\mathbb{R}$ of usual real numbers has the simple {\it self-similar} structure given by the order inverting bijective mapping \\

(4.1)~~~ $ \mathbb{R} \setminus ( -1, 1 ) \ni r \longmapsto \frac{1}{r} \in [ -1, 1 ] \setminus \{ 0 \} $ \\

This means that the {\it unbounded} set $( -\infty, -1 ] \cup [ 1, \infty )$ has the inverse order structure of the {\it bounded} set $[ -1, 0 ) \cup
( 0, 1]$, and of course, vice versa. Also (4.1) implies simply by {\it translation} the following self-similarity, centered not only around 0 as above,
but around every given $r_0 \in \mathbb{R}$, namely \\

(4.2)~~~ $ \mathbb{R} \setminus ( -1, 1 ) \ni r \longmapsto \frac{1}{r} + r_0 \in [ r_0 - 1, r_0 + 1 ] \setminus \{ r_0 \} $ \\

And further such self-similarities result by {\it scaling}, that is, by changing the unit, namely \\

(4.3)~~~ $ \mathbb{R} \setminus ( -u, u ) \ni r \longmapsto \frac{1}{r} + r_0 \in [ r_0 - \frac{1}{u}, r_0 + \frac{1}{u} ] \setminus \{ r_0 \} $ \\

for every $u > 0$. \\

Clearly, none of the above self-similarities refers to the structure at any given point $r_0 \in \mathbb{R}$, when such a point is considered itself
alone, but only to the structure of the sets \\

(4.4)~~~ $ [ r_0 - u, r_0 + u ] \setminus \{ r_0 \} = [ r_0 - u, r_0 ) \cup ( r_0, r_0 + u ],~~~ u > 0 $ \\

around such a point $r_0 \in \mathbb{R}$, sets which are whole neighbourhoods of $r_0$, from which, however, that point $r_0$ itself was eliminated. \\

This is certainly inevitable, since each point $r_0 \in \mathbb{R}$ is at a finite positive - thus non-infinitesimal - distance from any other point,
when considered in $\mathbb{R}$. Hence the moment one wants to focus on such a point $r_0 \in \mathbb{R}$ alone, one must exclude all other points in
$\mathbb{R}$, as there are {\it no} nonzero infinitesimals. \\

On the other hand, with the RPA-s, their self-similar structures are far more involved, due to the presence of their {\it infinitesimals}, and thus as
well, of their {\it infinitely large} elements. Indeed, this time, the self-similarities refer to each and every point itself {\it together} with its
infinitesimals, that is, the self-similarities refer to the whole monad of each such point. \\

Here, for simplicity, we shall indicate such self-similar structures in the particular case of reduced power fields, that is, RPF-s. Let therefore \\

(4.5)~~~ $ \mathbb{F}_{\cal U} = \mathbb{R}^\Lambda / {\cal I}_{\cal U} $ \\

be the RPF which corresponds to a given free ultrafilter ${\cal U}$ on $\Lambda$, and recall that such RPF-s are totally ordered. We also recall that
such RPF-s are strictly larger than $\mathbb{R}$. \\

Let us start with a self-similarity of RPF-s which does {\it not} exist in the case of the usual real line $\mathbb{R}$. Namely, it is easy to see that
we have the order inverting bijective mapping \\

(4.6)~~~ $ \mathbb{F}_{\cal U} \setminus Gal ( 0 ) \ni t \longmapsto \frac{1}{t} \in monad ( 0 ) \setminus \{ 0 \} $ \\

which means that the set of all {\it infinitely large} elements in $T_{\cal U}$ has the inverse order structure of the set of {\it infinitesimal}
elements from which one excludes 0. Also, through translation, we have, for each $t_0 \in T_{\cal U}$, the order inverting bijective mapping \\

(4.7)~~~ $ \mathbb{F}_{\cal U} \setminus Gal ( 0 ) \ni t \longmapsto \frac{1}{t} + t_0 \in monad ( t_0 ) \setminus \{ t_0 \} $ \\

which again are self-similarities {\it not} present in the case of the usual real line $\mathbb{R}$. \\

In the case of RPF-s, however, we have a far more rich possibility for scaling, since in addition to scaling with non-zero finite elements as in (4.3),
we can now also scale with all infinitely large elements, as well as with all infinitesimal elements, except for 0. \\

Let us consider the corresponding scalings, for a given $t_0 \in T_{\cal U}$. We take any $u \in \mathbb{F}_{\cal U},~ u > 0$ and obtain the order
inverting bijective mapping \\

(4.8)~~~ $ \mathbb{F}_{\cal U} \setminus ( - u, u ) \ni t \longmapsto
                            \frac{1}{t} + t_0 \in [ t_0 - \frac{1}{u}, t_0 + \frac{1}{u} ] \setminus \{ t_0 \} $ \\

where $\mathbb{F}_{\cal U} \setminus ( - u, u )$ will always contain {\it infinitely large} elements. \\

The difference with (4.3) is with respect to the sets \\

(4.9)~~~ $ [ t_0 - \frac{1}{u}, t_0 + \frac{1}{u} ] \setminus \{ t_0 \} $ \\

First of all, these sets are no longer mere subsets in $\mathbb{R}$, but instead, they are subsets in $\mathbb{F}_{\cal U}$, and will always contain
{\it infinitesimals}, since they contain nonvoid intervals. Furthermore, as seen below, they may on occasion also contain {\it infinitely large}
elements. \\

Also, $t_0$ and $u$ in (4.9) can independently be finite, infinitesimal, or infinitely large, thus resulting in 9 possible combinations and 6 distinct
outcomes regarding the set (4.9), which we list below. This is in sharp contradistinction with the case in (4.3) which applies to the real line
$\mathbb{R}$. Indeed : \\

1) Let us start listing the 9 different cases and 6 distinct outcomes with both $t_0$ and $u$ being {\it finite}. Then obviously (4.9) is a subset of
$Gal ( 0 )$, and it has the finite, non-infinitesimal length $\frac{2}{u}$. \\

2) When $t_0$ is {\it finite} and $u$ is {\it infinitesimal}, then the set (4.9) is infinitely large, and is no longer contained in $Gal ( 0 )$, however,
contains $Gal ( 0 ) \setminus \{ t_0 \}$. \\

3) If $t_0$ is {\it finite}, but $u$ is {\it infinitely large}, then (4.9) is again a subset of $Gal ( 0 )$, and in fact, it has the infinitesimal length
$\frac{2}{u}$, which means that it is a subset of $monad ( t_0 )$. \\

4) Let us now assume $t_0$ is {\it infinitesimal} and $u$ {\it finite}. Then regarding the set (4.9), we are back to case 1) above. \\

5) If both $t_0$ and $u$ are {\it infinitesimal} then the set (4.9) is as in 2) above. \\

6) When $t_0$ is {\it infinitesimal} and $u$ is {\it infinitely large}, the set (4.9) is as in 3) above. \\

7) Let us now take $t_0$ {\it infinitely large} and $u$ {\it finite}. Then the set (4.9) is disjoint from $Gal ( 0 )$, and it has the finite,
non-infinitesimal length $\frac{2}{u}$. \\

8) When $t_0$ {\it infinitely large} and $u$ {\it infinitesimal}, then the set (4.9) is again not contained in $Gal ( 0 )$, and it has the infinitely
large length $\frac{2}{u}$. Furthermore, depending on the relationship between $| t_0 |$ and $\frac{1}{u}$, it may, or it may not intersect $Gal ( 0
)$. \\

9) Finally, when both $t_0$ and $u$ are {\it infinitely large}, then the set (4.9) is disjoint from $Gal ( 0 )$, and it has the infinitesimal length
$\frac{2}{u}$. \\

Needless to say, the self-similar structure of RPA-s in general is still more rich and complex than the above in the case of RPF-s. \\ \\

{\bf 5. Conceiving Limits in Calculus} \\

Let us consider the fundamental operation in Calculus, namely, the limit of a sequence of real numbers \\

(5.1)~~~ $ \lim_{\,n \to \infty} x_n = x $ \\

where we assume that the respective $x \in \mathbb{R}$ cannot be defined by a finite information. Consequently, the terms $x_n$ in the above sequence may
have to contain unbounded information. \\

In this regard, what appears to be the essential novel phenomenon brought in by Calculus, when compared with the earlier Elementary Mathematics, is that
the left hand term in (5.1) means in a certain sense that

\begin{itemize}

\item infinitely many operations with unbounded amount of information are conceived in a finite usual time.

\end{itemize}

What is further remarkable in (5.1) is that in such a performance facilitated by Calculus there is {\it no} any kind of Zeno-type effect in usual time.
In fact, we may clearly note an opposite effect. Indeed, the terms in the sequence in (5.1) are not supposed to have a bounded, let alone, a decreasing
amount of information. Therefore, the effect of the above quite naturally is to ask the questions

\begin{itemize}

\item How can Calculus offer the possibility to conceive in finite usual time infinitely many operations with unbounded amount of information ?

\item Which is the kind of time structure within which Calculus manages to facilitate such a performance ?

\end{itemize}

And as if to further aggravate the situation, there comes the rigorous definition of the limit \\

(5.2)~~~ $ \lim_{\,n \to \infty} x_n = x $ \\

according to Calculus, namely \\

(5.3)~~~ $ \forall~ \epsilon > 0 ~:~ \exists~ m \in \mathbb{N} ~:~ \forall~ n \in \mathbb{N} ~:~ n \geq m
                                                                ~\Longrightarrow~ |\, x - x_n \,| \leq \epsilon $ \\

And the essential fact in (5.3) from the point of view of Calculus is that the two universal quantifiers $\forall$ which appear in it range over the
infinite domains $\epsilon \in ( 0, \infty )$ and $n \in \mathbb{N}$, respectively. Yet in the mind of a human being who knows and understands Calculus,
the operations of the respective two universal quantifiers happen in finite usual time. \\

But then, such a mental process in humans is obviously but a particular case of the ability of human mind to conceive infinity, be it actual or
potential, and do so in finite usual time. \\

One of the major novelties, therefore, brought about by Calculus, and specifically, by its quintessential operation of limit, is to place the {\it
pragmatic} aspect of infinity up front, and in fact, to highlight the human ability to deal with an actually infinite amount of arithmetical operations, and
do so in usual finite time. \\

In this regard, the ancient paradoxes of Zeno appear to be no more than an expression of a mental inability to make the very last step done by Newton,
namely, to jump from the potential infinity in a sequence $x_1, x_2, x_3, \ldots , x_n, \ldots $ to $\lim_{\,n \to \infty} x_n = x$, seen as an actual
infinity. \\

In this way, Calculus, within its specific realms, has given a first and major treatment of the age old issue of potential, versus actual infinity, and
has done so pragmatically, and in massively useful ways. \\

However, the most impressive theoretical approach to infinity have, so far, been the Set Theory of Cantor and Category Theory of Eilenberg and Mac
Lane. \\ \\

{\bf 6. Are Limits in Calculus Conceived in Monads of Time ?} \\

The above, and specifically the question

\begin{itemize}

\item Which is the kind of time structure within which Calculus manages to facilitate such a performance ?

\end{itemize}

may seem to lead as an answer to the conclusion that

\begin{itemize}

\item The process of conceiving a limit in Calculus must take place in a {\it monad} of time, thus in a non-Archimedean time structure.

\end{itemize}

One way to try to substantiate this conclusion is as follows. Let us take the limit process involved in (5.1) and assume a usual time interval $[ 0, T
]$, with a certain given $T > 0$, during which a human intelligence, competent in Calculus, grasps it. Certainly, we can conceive of an infinite
sequence, for instance \\

(6.1)~~~ $ 0 = t_0 < t_1 = T / 2 < t_2 = 2 T / 3 < \ldots < n T / ( n + 1 ) < \ldots $ \\

of moments of usual time which may correspond to the comprehension of the {\it presence} of the respective terms in (5.1), even if {\it not} necessarily
to their meaning, or in particular, to the amount of information they contain. After all, such a comprehension process of mere presence is nothing else
but what is supposed to be involved in the paradoxes of Zeno. \\

However, there seems to be an essential difference here. Namely, and as mentioned, the successive terms in (5.1) do not contain more a bounded amount of
information. Thus conceiving them not merely by their presence, it is not easy to see them accommodated within shorter and shorter usual time intervals,
such as for instance resulted from (6.1). And yet, this is precisely what Calculus does, or rather, what with the help of Calculus a competent human
intelligence can do. \\

And then, one possible explanation is that the respective unbounded amount of information in the successive terms in (5.1) is not being accommodated
either at the usual time moment $t_n$, or during the usual time interval $( t_n, t_{n+1} )$, but rather within the monad \\

(6.2)~~~ $ monad ( t_n ) $ \\

And needless to say, this may perfectly be possible, since in every RPF, each such monad, except for the single point $t_n$ itself, is a set of
uncountable cardinality which, according to (4.7), is self-similar with the set \\

(6.3)~~~ $ \mathbb{F}_{\cal U} \setminus Gal ( 0 ) $ \\

Consequently, each monad (6.2) can easily accommodate no matter how large an amount of finite information. \\ \\

{\bf 7. Is Space-Time Non-Archimedean ?} \\

To the extent that the above conclusion may be true with respect to the non-Archimedean structure of time, there is apparently no reason to further hold
to an Archimedean assumption related to space. \\

And then, we may as well conclude about a non-Archimedean structure for space-time. \\ \\

{\bf 8. Comments} \\

From the point of view of the Cartesian separation between "res extensa" and "res cogitans", the argument in section 6 above is questionable, since usual
time is supposed to belong to "res extensa", while the processes in human intelligence, and among them, those related to Calculus are assumed to be part
of "res cogitans". \\
Descartes is nowadays widely accused for being a dualist, forgetting completely that, as so many prominent scientists of his age, he was deeply religious.
And as such, he could of course in no way be dualist, since he saw God as the underlying sole source and support of all Creation, including, no doubt,
both "res extensa" and "res cogitans". \\

A great merit, however, and one that is still not quite realized, let alone accepted today, in that Cartesian separation of realms is that the thinking
human intelligence is not automatically taken out of all consideration, with its focus being reduced exclusively to one or another aspect of "res
extensa". Indeed, "res cogitans" - much unlike "res extensa" - must by its own very definition relate to itself as well in its essential function of
thinking, and not only to "res extensa", this self-referentiality being thus in its nature. \\
Nowadays however, in spite of the fact that for more than eight decades by now the Copenhagen Interpretation of Quantum Theory keeps warning of the
inadequacy of considering the observer and the observation as being totally outside of, and not directly affecting the quantum physical process focused
upon, the general scientific practice is still to separate - and then disregard - totally the issues related to human intelligence, while it is involved
in scientific thinking. \\

And to give a simple example from modern, post-Calculus Mathematics, we can ask the following question

\begin{itemize}

\item How does it happen that human intelligence can conceive in finite usual time of sets of immensely larger cardinal than that of a usual time
interval ?

\end{itemize}

Needless to say, there are any number of such or similar questions, some of them were mentioned in [3]. \\

A possible demerit of the Cartesian approach is in the assumed sharp qualitative difference between the mentioned two realms. In this regard it is worth
noting that Quantum Theory, for instance, in whichever of its many interpretations tends to agree that it is not so easy to define precisely where is the
separation between a quantum process and the measuring apparatus employed in its observation. \\

As for the argument in section 6 above, the more one distances oneself from any Cartesian type separation of realms, the more that argument may gain in
strength. \\

By the way of quanta, one may also remark that not only the RPF-s, but the more general and yet more rich and complex RPA-s themselves may possibly be
conceived as modelling space-time, even if they would lead to infinite dimensions not only for space, but also for time. And in such models, precisely
because of the considerably increased richness and complexity in the available structure, one may eventually find a possibility to place the Many-Worlds
interpretation of Quantum Theory according to Everett.

\newpage

{\LARGE \bf Part 2} \\ \\

{\LARGE \bf Walkable Worlds : Universes next \\ \\ to and/or within Universes ... and \\ \\ so on ad infinitum ...} \\ \\

{\bf Abstract} \\

There is an insufficient awareness about the {\it rich and complex} structure of various totally ordered scalar fields
obtained through the {\it ultrapower} construction. This rich and complex structure comes from the presence of {\it
infinitesimals} in such fields, presence which leads to the fact that such fields are {\it non-Archimedean}. Here, with
the concept of {\it walkable world}, which has highly intuitive and pragmatic geometric meaning, the mentioned rich and
complex structure is illustrated. The issues presented have relevance for what are usually called the "infinities in
physics". \\ \\

{\bf 1. Totally Ordered Fields as Ultrapowers} \\

There is a simple way to construct totally ordered fields which contain the usual field given by the real line
$\mathbb{R}$ as a rather small and subset. This construction is a particular case of what is well known in Model Theory, a
branch of Mathematical Logic, as {\it reduced powers}, and proceeds as follows,[2-8]. Given an infinite index set
$\Lambda$, we take on it any {\it ultrafilter} ${\cal U}$ which satisfies the condition \\

(1.1)~~~ $ {\cal F}_{re} ( \Lambda ) \subseteq {\cal U} $ \\

where \\

(1.2)~~~ $ {\cal F}_{re} ( \Lambda ) = \{~ I \subseteq \Lambda ~~|~~ \Lambda \setminus I ~~\mbox{is finite} ~\} $ \\

is called the Frech\`{e}t filter on $\Lambda$. Further, we define on $\mathbb{R}^\Lambda$ the equivalence relation
$\approx_{\cal U}$ by \\

(1.3)~~~ $ x \approx_{\cal U} y ~~\Longleftrightarrow~~
                        \{~ \lambda \in \Lambda ~|~ x ( \lambda ) = y ( \lambda ) ~\} \in {\cal U} $ \\

Finally, through the usual quotient construction, we obtain the {\it ultrapower field} \\

(1.4)~~~ $ \mathbb{F}_{\cal U} = \mathbb{R}^\Lambda / \approx_{\cal U} $ \\

which proves to have the following two properties. The mapping \\

(1.5)~~~ $ \mathbb{R} \ni r \longmapsto ( u_r )_{\cal U} \in \mathbb{F}_{\cal U} $ \\

is an {\it embedding of fields} in which $\mathbb{R}$ is a strict subset of $\mathbb{F}_{\cal U}$, where $u_r \in
\mathbb{R}^\Lambda$ is defined by $u_r ( \lambda ) = r$, for $\lambda \in \Lambda$, while $( u_r )_{\cal U}$ is the coset
of $u_r$ with respect to the equivalence relation $\approx_{\cal U}$. For simplicity we shall denote $( u_r )_{\cal U} =
r$, for $r \in \mathbb{R}$, and thus (1.5) takes the form \\

(1.6)~~~ $ \mathbb{R} \ni r \longmapsto r \in \mathbb{F}_{\cal U},~~~
                         \mbox{or simply}~~~ \mathbb{R} \subsetneqq \mathbb{F}_{\cal U} $ \\

Further, on $\mathbb{F}_{\cal U}$ we have the {\it total order} \\

(1.7)~~~ $ ( x )_{\cal U} \leq ( y )_{\cal U} ~~\Longleftrightarrow~~
             \{~ \lambda \in \Lambda ~|~ x ( \lambda ) \leq y ( \lambda ) ~\} \in {\cal U} $ \\

where $x, y \in \mathbb{R}^\Lambda$. \\

It should be noted that the nonstandard reals $^*\mathbb{R}$ are a particular case of the above ultrapower fields
(1.4). \\

The general case of the above construction is that of {\it reduced power algebras}, which goes as follows. Let ${\cal F}$
be any filter on $\Lambda$ which satisfies \\

(1.8)~~~ $ {\cal F}_{re} ( \Lambda ) \subseteq {\cal F} $ \\

We define on $\mathbb{R}^\Lambda$ the corresponding equivalence relation $\approx_{\cal F}$ by \\

(1.9)~~~ $ x \approx_{\cal U} y ~~\Longleftrightarrow~~
                        \{~ \lambda \in \Lambda ~|~ x ( \lambda ) = y ( \lambda ) ~\} \in {\cal F} $ \\

Then, through the usual quotient construction, we obtain the {\it reduced power algebra} \\

(1.10)~~~ $ \mathbb{A}_{\cal F} = \mathbb{R}^\Lambda / \approx_{\cal F} $ \\

which has the following two properties. The mapping \\

(1.11)~~~ $ \mathbb{R} \ni r \longmapsto ( u_r )_{\cal F} \in \mathbb{A}_{\cal F} $ \\

is an {\it embedding of algebras} in which $\mathbb{R}$ is a strict subset of $\mathbb{A}_{\cal F}$, where $u_r \in
\mathbb{R}^\Lambda$ is defined by $u_r ( \lambda ) = r$, for $\lambda \in \Lambda$, while $( u_r )_{\cal F}$ is the coset
of $u_r$ with respect to the equivalence relation $\approx_{\cal F}$. Further, on $\mathbb{A}_{\cal F}$ we have the {\it
partial order} \\

(1.12)~~~ $ ( x )_{\cal F} \leq ( y )_{\cal F} ~~\Longleftrightarrow~~
             \{~ \lambda \in \Lambda ~|~ x ( \lambda ) \leq y ( \lambda ) ~\} \in {\cal F} $ \\

where $x, y \in \mathbb{R}^\Lambda$. \\ \\

{\bf 2. Walkable Worlds ...} \\

An essential property of the ultrapower fields (1.4) is that they are {\it no} longer Archimedean. In other words, unlike
the usual field given by the real line $\mathbb{R}$, they do {\it not} satisfy the Archimedean Axiom \\

(2.1)~~~ $ \exists~~ u > 0 ~:~ \forall~~ v > 0 ~:~ \exists~~ n \in \mathbb{N} ~:~ n u > v $ \\

As it happens, and it is still seldom realized, the fact that the ultrapower fields (1.4) are {\it non-Archimedean}, and
also are {\it larger} than $\mathbb{R}$, gives them an extremely rich and complex both {\it local} and {\it global}
structure. Here we shall illustrate that fact with the help of the concept of {\it walkable world}, a concept which is
highly intuitive in its pragmatic and geometric meaning. \\

Let us briefly recall some of the basic features of the ultrapower fields (1.4). One that follows immediately from the
fact that they are non-Archimedean is that their elements $t \in \mathbb{F}_{\cal U}$ are of {\it three} kind, namely,
{\it infinitesimal, finite}, and {\it infinitely large}, defined by the following respective conditions \\

(2.2)~~~ $ \forall~~ r \in \mathbb{R},~ r > 0 ~:~ t \in ( - r, r ) $ \\

(2.3)~~~ $ \exists~~ r \in \mathbb{R},~ r > 0 ~:~ t \in ( - r, r ) $ \\

(2.4)~~~ $ \forall~~ r \in \mathbb{R},~ r > 0 ~:~ t \notin ( - r,  r ) $ \\

where for $a, b \in \mathbb{F}_{\cal U}$, we denote $( a, b ) = \{ s \in \mathbb{F}_{\cal U} ~|~ a < s < b \}$. Now,
following Leibniz, one denotes \\

(2.5)~~~ $ monad ( 0 ) = \{~ t \in \mathbb{F}_{\cal U} ~|~ t ~\mbox{infinitesimal}~ \} $ \\

and calls it the {\it monad} of $0 \in \mathbb{F}_{\cal U}$, while following Keisler, [1], one denotes \\

(2.6)~~~ $ Gal ( 0 ) = \{~ t \in \mathbb{F}_{\cal U} ~|~ t ~\mbox{finite}~ \} $ \\

and calls it the {\it Galaxy} of $0 \in \mathbb{F}_{\cal U}$. It follows that \\

(2.7)~~~ $ Gal ( 0 ) = \bigcup_{r \in \mathbb{R}} monad ( r ) $ \\

where for $t \in \mathbb{F}_{\cal U}$, we denote \\

(2.8)~~~ $ monad ( t ) = t + monad ( 0 ) $ \\

Finally \\

(2.9)~~~ $ \mathbb{F}_{\cal U} \setminus Gal ( 0 ) $ \\

is the set of infinitely large elements in the ultrapower field $\mathbb{F}_{\cal U}$. \\

And now, to the walkable worlds ... \\

Let $t, u \in \mathbb{F}_{\cal U}$, with $u > 0$. Then we denote \\

(2.10)~~~ $ \begin{array}{l}
                WW ( t, u ) =
                  \{~ s \in \mathbb{F}_{\cal U} ~|~ \exists~ n \in \mathbb{N} ~:~ s \in ( t - n u, t + n u ) ~\} = \\ \\
                  ~~~~~~~~~~~~~ = \bigcup_{\, n \in \mathbb{N}}~ ( t - n u, t + n u )
            \end{array} $ \\

which is the set of elements $s \in \mathbb{F}_{\cal U}$ that can be reached, starting at $t$, by a {\it finite} number of
steps of length $u$. Thus $WW ( t, u )$ is the {\it walkable world} around $t$, by steps of length $u$. \\

Obviously \\

(2.11)~~~ $ Gal ( 0 ) = WW ( 0, 1 ) $ \\

and \\

(2.12)~~~ $ s \in WW ( t, u ) ~~\Longrightarrow~~ WW ( s, u ) = WW ( t, u ) $ \\

while \\

(2.13)~~~ $ WW ( t, v ) = WW ( t, u ) $ \\

for all $v \in \mathbb{F}_{\cal U},~ v > 0$, such that \\

(2.14)~~~ $ m u \leq v \leq ( m + 1 ) u $ \\

for some $m \in \mathbb{N}$. \\

Furthermore, for every $t, u \in \mathbb{F}_{\cal U}$, with $u > 0$, we have the {\it order isomorphism} \\

(2.15)~~~ $ WW ( t, u ) \ni s \longmapsto ( s - t ) / u \in WW ( 0, 1 ) $ \\

thus any two walkable worlds are {\it order isomorphic}. \\

The non-Archimedean nature of $\mathbb{F}_{\cal U}$ results, among others, in \\

(2.16)~~~ $ monad ( 0 ) \neq WW ( t, u ),~~~ \mathbb{F}_{\cal U} \neq WW ( t, u ) $ \\

and in fact \\

(2.17)~~~ $ WW ( t, u ) $ is an infinitely small part of $ \mathbb{F}_{\cal U} $ \\

for all $t, u \in \mathbb{F}_{\cal U}$, with $u > 0$. \\

{\bf Remark 2.1.} \\

In view of (2.11), in any given ultrapower field $\mathbb{F}_{\cal U}$, the whole of $Gal ( 0 )$ is only one single
walkable world, while as seen in the next section, there are infinitely many walkable worlds, either disjoint from one
another, or nested within one another. \\
Here is already one of the essential differences with ultrapower fields, when compared with the usual field of the real
line $\mathbb{R}$. Indeed, in the latter case, what corresponds to $Gal ( 0 ) = WW ( 0, 1 )$ is the whole of $\mathbb{R}$,
namely, for every $t, u \in \mathbb{R},~ u > 0$, we have \\

(2.18)~~~ $ \{~ s \in \mathbb{R} ~|~ \exists~ n \in \mathbb{N} ~:~ s \in ( t - n u, t + n u ) ~\} = \mathbb{R} $ \\

In other words, the usual field of the real line $\mathbb{R}$ is but only one single walkable world, while in the case of
ultrapower fields, each walkable world is merely an infinitely small part of such a field, as seen in (2.17). \\

So much limitation is imposed upon the structure of $\mathbb{R}$ by the acceptance of the Archimedean Axiom. \\ \\

{\bf 3. Universes Next To and/or Within Universes, and so on \\
        \hspace*{0.45cm} Ad Infinitum ...} \\

Let us see now the way two arbitrary walkable worlds $WW ( t, u )$ and $WW ( s, v )$, with $t, u, s, v \in
\mathbb{F}_{\cal U},~ u, v > 0$, can relate to one another. As it turns out, we can distinguish the following relative
situations : \\

(3.1)~~~ $ WW ( t, u ) = WW ( s, v ) $ \\

(3.2)~~~ $ WW ( t, u ) \cap WW ( s, v ) = \phi $ \\

(3.3)~~~ $ WW ( t, u ) \cap WW ( s, v ) \neq \phi,~~~ WW ( t, u ) \neq WW ( s, v ) $ \\

and in the last case, we have either the nesting \\

(3.3.1)~~~ $ WW ( t, u ) $ is an infinitesimal part of $ WW ( s, v ) $ \\

or the nesting \\

(3.3.2)~~~ $ WW ( s, v ) $ is an infinitesimal part of $ WW ( t, u ) $ \\

Furthermore, concerning (3.2), there are infinitely many walkable worlds which are pair-wise disjoint. As for (3.3), the
respective nestings in (3.3.1) and (3.3.2) go on in each case infinitely many times. \\

Let us go more in detail on the possible relation between two walkable worlds in the case (3.3). We note that, given $u, v
\in \mathbb{F}_{\cal U},~ u, v > 0$, we can in view of the total order on $\mathbb{F}_{\cal U}$ always assume that $ v
\leq u$. Thus the following two alternatives result \\

(3.4)~~~ either $u / v$ is finite, or $u / v$ is infinitely large \\

and correspondingly, either \\

(3.5)~~~ $ \exists~~ n \in \mathbb {N},~ n \geq 1 :~ u \leq n v $ \\

or \\

(3.6)~~~ $ \forall~~ n \in \mathbb {N},~ n \geq 1 :~ n v \leq u $ \\

Now in view of (2.13), (2.14), the alternative (3.5) is not compatible with (3.3), since it leads to (3.1). Thus we remain
with (3.6). And then (3.3.2) follows. \\
What is important to note here is that $u$ and $v$ can be in the following 6 situations \\

1) both $u$ and $v$ are infinitesimal \\

2) $u$ is finite, $v$ is infinitesimal \\

3) both $u$ and $v$ are finite \\

4) $u$ is infinitely large, $v$ is infinitesimal \\

5) $u$ is infinitely large, $v$ is finite \\

6) both $u$ and $v$ are infinitely large \\ \\

{\bf 4. Infinities in Physics ?} \\

In sharp contrast with the rich and complex structure of ultrapower fields, which by necessity are {\it not} Archimedean,
and thus exhibit the mentioned wealth of {\it walkable worlds}, all the known theories of Physics are, and have been stuck
into only one {\it single} such walkable world. No wonder that so called "infinities in physics" trouble various theories
of Physics, and furthermore, lead to highly questionable ad-hock attempts at solution. \\

On the other hand, whatever proves to be "infinity" in Physics becomes just another usual and regular element in the
ultrapower fields, thus all algebraic operations can be effectuated upon such an element, without any concern or
restrictions. \\

And then, all that takes is simply to discard the Archimedean Axiom. \\

\end{document}